\documentclass[11pt]{article}
\usepackage{amsfonts}
\usepackage{mathrsfs}
\usepackage{amsmath}
\usepackage{amssymb}
\usepackage{graphicx}
\usepackage{epsfig}
\usepackage{epic}
\usepackage{cite}
\usepackage{setspace}
\usepackage{amsthm,latexsym}
\usepackage{float}

\renewcommand{\paragraph}{\roman{paragraph}}
\newtheorem{theorem}{\scshape \mdseries  Theorem}[section]
\newtheorem{lemma}[theorem]{\scshape \mdseries  Lemma}
\newtheorem{coro}[theorem]{\scshape \mdseries  Corollary}

\newtheorem{prop}[theorem]{\scshape \mdseries  Proposition}

\begin{document}

\title{\sf Automorphism group of the subspace\\ inclusion graph of a vector space\thanks{}}
\author{Dein Wong \thanks{Corresponding
author. E-mail address:wongdein@163.com. Supported by ``the National Natural Science Foundation of China (No.11571360)".},\ \ \ \ Xinlei Wang,\ \ \ \  Fenglei Tian\\
{\small  \it  School of Mathematics, China University of Mining and
Technology, Xuzhou 221116,  China.}
   }
\date{}
\maketitle

\noindent {\bf Abstract:}\ \
 In a recent  paper [Comm. Algebra, 44(2016) 4724-4731],  Das  introduced the graph $\mathcal{I}n(\mathbb{V})$, called subspace inclusion graph on a finite dimensional vector space $\mathbb{V}$, where the vertex set is the collection of nontrivial proper subspaces of $\mathbb{V}$ and two vertices are adjacent if one is properly contained in another.   Das studied the diameter, girth, clique number, and chromatic number of  $\mathcal{I}n(\mathbb{V})$ when the base field is arbitrary, and he also   studied   some other properties of  $\mathcal{I}n(\mathbb{V})$  when the base field is finite. In this paper, the automorphisms of  $\mathcal{I}n(\mathbb{V})$ are determined when the base field is finite.



\section{Introduction}
\quad
Let $\mathbb{V}$ be a finite dimensional vector space over a field $\mathbb{F}$ of dimension greater than $1$.  Das in [5] defined  the  subspace inclusion graph $\mathcal{I}n(\mathbb{V})$ of $\mathbb{V}$  as follows: The vertex set $V$ of $\mathcal{I}n(\mathbb{V})$ is the collection of nontrivial proper  subspaces of $\mathbb{V}$ and for $W_1, W_2 \in V$, $W_1$ is adjacent to $W_2$, written as $W_1\sim W_2$,  if either $W_1 \subset W_2$ or $W_2\subset W_1$. The main results of [5] are as follows:

\begin{prop} (1) The diameter of  $\mathcal{I}n(\mathbb{V})$ is $3$ if the dimension of $\mathbb{V}$ is at least $3$;\\
(2) The girth of  $\mathcal{I}n(\mathbb{V})$ is either $3$, $6$ or $\infty$;\\
(3) The clique number and the chromatic number of   $\mathcal{I}n(\mathbb{V})$  are both $dim(\mathbb{V})-1$.
\end{prop}

If the base field $\mathbb{F}$ is a finite field with $q$ elements, the author obtained the following result  about the vertex degrees of  $\mathcal{I}n(\mathbb{V})$, which will be applied in our  result.

\begin{lemma} {\rm (Theorem 6.1, [5])}   If $W$ is a $k$-dimensional nontrivial proper subspace of $\mathbb{V}$, then the degree of $W$ is  $$\sum_{i=1}^{k-1}\left[\begin{array}{c}k\\i\end{array}\right]_q+\sum_{i=1}^{n-k-1}\left[\begin{array}{c}n-k\\i\end{array}\right]_q,$$
where $$\left[\begin{array}{c}n\\k\end{array}\right]_q=\frac{(q^n-1)(q^{n-1}-1)\cdots(q^{n-k+1}-1)}{(q^k-1)(q^{k-1}-1)\cdots(q-1)}.$$
\end{lemma}

Till date, a lot of researches, e.g., [1,3,4,5,7] have been performed in connecting graph structures to   subspaces of vector spaces.

Automorphisms  of graphs are important in algebraic graph theory because they reveal the relationship  between the vertices of the graph. Automorphisms of the zero-divisor graph of a ring  $R$, denoted by $\Gamma(R)$ and defined as a graph with $Z(R)\setminus \{0\}$ as vertex set and there is a directed edge from a vertex $x$ to a distinct vertex $y$ if and only if $xy=0$, have attracted a lot of attention. In [2], Anderson and Livingston proved
that the automorphism group of the zero-divisor graph of $\mathbb{Z}_n$ is a  direct product of symmetric groups for $n \geq 4$ a
nonprime integer. It was shown in [6] that Aut$(\Gamma(R)$) is isomorphic to the symmetric
group of degree $p + 1$, when $R =$ $M_2(\mathbb{Z}_p)$ with $p$ a prime.  Park and  Han [8]
 proved that Aut$(\Gamma(R))\cong  S_{q+1}$ for $R = $ $M_2(F_q)$ with $F_q$ an arbitrary  finite field.
In  [11],  Wong et al. determined the automorphisms of the zero-divisor graph with vertex set of all rank one upper triangular matrices over a finite field. By applied the main theorem in [11], Wang [9] and [10] respectively determined the automorphisms of the zero-divisor graph defined on all $n\times n$ upper triangular matrices  or on all $n\times n$ full matrices when the base field is finite.

Now, a natural problem arises:{\sf What about the automorphisms of $\mathcal{I}n(\mathbb{V})$?}

Denote by $dim(W)$ the dimension of a subspace $W$ of $\mathbb{V}$.
If   $dim(\mathbb{V})=1$, then the vertex set $V$ of  $\mathcal{I}n(\mathbb{V})$ is empty; if  $dim(\mathbb{V})=2$ then   $\mathcal{I}n(\mathbb{V})$ is a graph consisting of some isolated vertices, thus any permutation on $V$ is an automorphism of $\mathcal{I}n(\mathbb{V})$. If $dim(\mathbb{V})\geq 3$, the situation is quite different. One will find that some nontrivial automorphisms do exist. In this paper, we solve the above problem for the case when  the base field is finite; if the base field if infinite, the problem is left open.

If no other explicit  mention,  $\mathcal{I}n(\mathbb{V})$ is an $n$-dimensional vector space over a finite field $F_q$ of $q$ elements, where $q=p^m$ and $p$ is a prime integer. Let $\{\epsilon_1, \epsilon_2, \ldots, \epsilon_n\}$  be a   set of base  of   $\mathcal{I}n(\mathbb{V})$. Every $\alpha\in \mathbb{V}$ can be uniquely written as $\alpha=\sum_{i=1}^na_i\epsilon_i$ with $a_i\in F_q$. If $\alpha=\sum_{i=1}^na_i\epsilon_i, \beta=\sum_{i=1}^nb_i\epsilon_i$ satisfy $\sum_{i=1}^n a_ib_i=0$, we write $\alpha\perp \beta$ to denote that  $\alpha, \beta$ are orthogonal. For a subspace $W$ of $\mathbb{V}$, set $$W^\perp=\{\alpha\in \mathbb{V}:\alpha\perp w,\ \forall w\in W\}.$$ It is easy to see that:\\
(i)\  $W^\perp$ is a subspace of $\mathbb{V}$; \ \ \
(ii) \ $(W^\perp)^\perp=W$;\ \ \
(iii)\  $dim(W)+dim(W^\perp)=dim(\mathbb{V})$.\\
We reminder the reader that $\mathbb{V}$ possibly fails to be the sum of $W$ and $W^\perp$ since the intersection of $W$ and $W^\perp$ possibly contains nonzero vectors.
Before announcing the main result of this paper, we introduce three standard automorphisms of  $\mathcal{I}n(\mathbb{V})$.

\vskip 1.5mm
  {\bf Involution of   $\mathcal{I}n(\mathbb{V})$}

Let  $\tau$ be the mapping from $V$, the vertex set of  $\mathcal{I}n(\mathbb{V})$, to itself, which sends any $W\in V$ to $W^\perp$. Then $\tau^2$ is the identity mapping on ${V}$, which implies that $\tau$ is a bijection on $V$.
Noting that  $W_1\subset W_2$ if and only if  $W_2^\perp\subset W_1^\perp$, we have $W_1\sim W_2$ if and only if $\tau(W_1)\sim \tau(W_2)$. Consequently,  $\tau$ is an automorphism of $\mathcal{I}n(\mathbb{V})$, which is called the involution of  $\mathcal{I}n(\mathbb{V})$.

\vskip 1.5mm
  {\bf Invertible linear transformation}

Let $X=[x_{ij}]$ be an $n\times n$ invertible matrix over $F_q$.  Then the mapping $\theta_X$ on $\mathbb{V}$ defined by $$\theta_X(\sum_{i=1}^n b_i\epsilon_i)=\sum_{i=1}^n (\sum_{j=1}^nx_{ij}b_j)\epsilon_i,\ \forall\ \sum_{i=1}^n b_i\epsilon_i\in \mathbb{V},$$
 is an invertible linear transformation on $\mathbb{V}$. The mapping from $V$ to itself,  also written as $\theta_X$,  sending any $W\in V$ to $ \{\theta_X(w):w\in W\}$ is an automorphism of $\mathcal{I}n(\mathbb{V})$.

\vskip 1.5mm
  {\bf Field automorphism}

Let $f$ be an automorphism of the base field $F_q$ and let  $\chi_f$ be the mapping on $\mathbb{V}$ defined by $$\chi_f(\sum_{i=1}^n b_i\epsilon_i)=\sum_{i=1}^nf(b_i)\epsilon_i,\ \forall\ \sum_{i=1}^nb_i\epsilon_i\in \mathbb{V}.$$
Then the related mapping, also written as $\chi_f$,  on $V$ sending any $W\in V$ to $\{\chi_f(w):w\in W\}$ is an automorphism of $\mathcal{I}n(\mathbb{V})$, which is called a field automorphism of $\mathcal{I}n(\mathbb{V})$.

   The main result of this paper is the following theorem.
   \begin{theorem} If $F$ is a finite field and  $dim(\mathbb{V})\geq 3$, then a mapping  $\sigma$ on $V$ is an automorphism of $\mathcal{I}n(\mathbb{V})$ if and only if $\sigma$ can be uniquely decomposed  either as $\sigma= \theta_X\circ\chi_f $ or as $\sigma= \tau\circ\theta_X\circ\chi_f$, where $\tau, \theta_X$ and $\chi_f$ are as just defined.
\end{theorem}

\begin{coro} Let $F=F_q$ with $q=p^m$, where $p$ is a prime integer. If $dim(\mathbb{V})\geq 3$,  then the automorphism group of $\mathcal{I}n(\mathbb{V})$ is isomorphic to $\mathbb{Z}_2\times PGL_n(F_q)\times \mathbb{Z}_m$ , where $PGL_n(F_q)$ is the quotient group of all $n\times n$ invertible matrices over $F_q$ to the normal subgroup of all nonzero scalar matrices over $F_q$.
\end{coro}

We will give a proof for   theorem 1.3 in the next section.

\section{ Proof of Theorem 1.3 }

\quad The subspace spanned by a subset $S$ of $\mathbb{V}$ is denoted by $[S]$.
An $1$-dimensional subspace $W$ can be written as  $W=[\alpha]$ with $\alpha=\sum_{i=1}^na_n\epsilon$ a nonzero vector of $\mathbb{V}$. As $[\alpha]=[b\alpha ]$ for any $0\not=b\in F_q$, the expression  is not unique. However, if we require the first nonzero coefficient of $\alpha$ is $1$, that is  $\alpha$ is of the form $$\alpha=\epsilon_k+a_{k+1}\epsilon_{k+1}+\ldots+a_n\epsilon_n,\ \  1\leq k\leq n,$$ then  the expression $W=[\alpha]$ is unique.
Such an expression of an $1$-dimensional subspace  is called standard. In what follows, all $1$-dimensional subspaces will be expressed in standard form. The degree of a vertex $W\in V$ in  $\mathcal{I}n(\mathbb{V})$ is denoted by $d(W)$. Before giving a proof for Theorem 1.3 we introduce some lemmas for latter use. The first  one follows from Lemma 1.2 immediately.

\begin{lemma}  For   $W, W'\in V$, we have $d(W)=d(W^\perp)$, and $d(W')<d(W)$ if  $1\leq dim(W)<dim(W')\leq \frac{n}{2}$.
\end{lemma}

\noindent{\sf Proof.}\ Suppose that $dim(W)=k$. Then $dim(W^\perp)=n-k$, thus by Corollary 6.2 in [5], we have $d(W)=d(W^\perp)$.
If $2\leq k\leq \frac{n}{2}$, then it follows from Lemma 1.2 that the degree of a $(k-1)$-dimensional subspace $U$ and the degree of a $k$-dimensional subspace $U'$ respectively are  $$d(U)=\sum_{i=1}^{k-2}\left[\begin{array}{c}k-1\\i\end{array}\right]_q+\sum_{i=1}^{n-k}\left[\begin{array}{c}n-k+1\\i\end{array}\right]_q,$$ $$d(U')=\sum_{i=1}^{k-1}\left[\begin{array}{c}k\\i\end{array}\right]_q+\sum_{i=1}^{n-k-1}\left[\begin{array}{c}n-k\\i\end{array}\right]_q.$$
By direct computation  we find that    $d(U)>d(U')$, from which  the second assertion follows.  \hfill$\square$

\begin{lemma}  Let $\sigma$ be an automorphism of $\mathcal{I}n(\mathbb{V})$. If $\sigma$ sends an $1$-dimensional subspace    to an $1$-dimensional subspace, then $\sigma$  sends every $k$-dimensional subspace to a subspace of equal dimension, where $1\leq k\leq n-1$.
\end{lemma}

\noindent{\sf Proof.} Suppose $\sigma$ sends  an $1$-dimensional subspace  $[\alpha_1]$ to an   $1$-dimensional subspace.   Firstly, we consider the  case when $k=1$.
Let $[\alpha]$ be an   $1$-dimensional subspace different from $[\alpha_1]$. Expand $\alpha_1$ to a set of base  of $\mathbb{V}$ as $\alpha_1, \alpha, \alpha_3, \ldots, \alpha_n$, where $n=dim(\mathbb{V})$. Then the following $n-1$ vertices:
 $$[\alpha_1], \ [\alpha_1,\alpha],\  [\alpha_1,\alpha,\alpha_3], \ldots, [\alpha_1,\alpha,\alpha_3,\ldots,\alpha_{n-1}]$$
  forms a maximum clique of  $\mathcal{I}n(\mathbb{V})$ (see the proof for Theorem 5.1 in [5]). The image of this clique under $\sigma$ is also a  maximum clique of  $\mathcal{I}n(\mathbb{V})$. Since $\sigma$ sends $[\alpha_1]$ to an $1$-dimensional subspace of $\mathbb{V}$, the image of $[\alpha_1,\alpha,\alpha_3,\ldots,\alpha_{n-1}]$ under $\sigma$ must   be an $(n-1)$-dimensional subspace of $\mathbb{V}$ (otherwise, it follows from Lemma 2.1 that  the dimension of this space is $1$ and thus this space is not adjacent to the image of $[\alpha_1]$ under $\sigma$, a contradiction). The following $n-1$ vertices:
 $$[\alpha],\  [\alpha_1,\alpha],\  [\alpha_1,\alpha,\alpha_3], \ldots, [\alpha_1,\alpha,\alpha_3,\ldots,\alpha_{n-1}]$$
 also induce   a maximum clique of  $\mathcal{I}n(\mathbb{V})$. Since the image of this clique under $\sigma$ is also a  maximum clique of  $\mathcal{I}n(\mathbb{V})$ and  the dimension of the image of $[\alpha_1,\alpha,\alpha_3,\ldots,\alpha_{n-1}]$ under $\sigma$ is  $(n-1)$, we confirm that the  image of $[\alpha]$ under $\sigma$ is of $1$-dimensional (otherwise, such a space is of $(n-1)$-dimensional and thus it is not adjacent to the image of  $[\alpha_1,\alpha,\alpha_3,\ldots,\alpha_{n-1}]$, a contradiction).

Next, We proceed by induction on $k$, for $1\leq k\leq n-1$,  to prove that $\sigma$ sends any $k$-dimensional subspace to a subspace of equal dimension.
The case when $k=1$ has been proved.   Suppose $\sigma$ sends every $i$-dimensional subspace to a subspace of equal dimension for $1\leq i\leq k-1$.
Let $W$ be a $k$-dimensional subspace. We confirm the dimension of   $ \sigma(W)$ is at least $k$.  Otherwise   the dimension of $\sigma(W)$ equals to the dimension of a proper subspace, say $U$, of $W$, thus
$\sigma(W)$ and $\sigma(U)$ are not adjacent since they have the same dimension (the induction hypothesis implies that $dim(\sigma(U))=dim(U)$), which is a contradiction to $U\sim W$. Hence, $dim(\sigma(W))\geq dim(W)$. Observing that $\sigma^{-1}$ also sends an   $1$-dimensional subspace    to an $1$-dimensional subspace and it sends $\sigma(W)$ to $W$, we have
$$dim(W)=dim(\sigma^{-1}(\sigma(W)))\geq dim(\sigma(W)).$$
Consequently,  $dim(\sigma(W))=dim(W)$. \hfill$\square$

\begin{lemma}  Let $\sigma$ be an automorphism of $\mathcal{I}n(\mathbb{V})$ which sends an $1$-dimensional subspace of $\mathcal{I}n(\mathbb{V})$   to an $1$-dimensional subspace, $W$   a $k$-dimensional subspace of $\mathbb{V}$ with a set of  base $\alpha_1, \alpha_2, \ldots, \alpha_k$ and suppose $\sigma([\alpha_i])=[\beta_i]$ for  $1\leq i\leq k$. If $\beta_1, \beta_2, \ldots, \beta_k$ are  linearly independent, then $\sigma(W)=[\beta_1, \beta_2, \ldots, \beta_k]$.
\end{lemma}

\noindent{\sf Proof.}   Assume that $k\geq 2$. For $1\leq i\leq k$, applying $\sigma$ to  $[\alpha_i]\subset W$ we have $[\beta_i]\subset \sigma(W)$, which implies that $[ \beta_1, \beta_2, \ldots, \beta_k]\subseteq \sigma(W)$. By Lemma 2.2, the dimension of $\sigma(W)$ is $k$, which  equals  the dimension of $[ \beta_1, \beta_2, \ldots, \beta_k]$, thus we have   $\sigma(W)=[\beta_1, \beta_2, \ldots, \beta_k]$.\hfill$\square$

Let $I$ be the $n\times n$ identity matrix and let  $E_{ij}$, for $1\leq i,j\leq n$, be the $n\times n$ matrix unit  whose $(i,j)$ position is $1$ and all other positions are $0$.

\begin{lemma}  Let $\sigma$ be an automorphism of $\mathcal{I}n(\mathbb{V})$, $1\leq k\leq n$. If $\sigma$ respectively fixes every  $[\epsilon_i]$ for $i=1,2, \ldots, k-1$, then  there exists an invertible matrix $A $ such that  $\theta_{A }\circ\sigma$  fixes  every  $[\epsilon_j]$ for $j=1,2, \ldots, k$.
\end{lemma}

\noindent{\sf Proof.} We only consider the case when $2\leq k$, the case when $k=1$ can be proved similarly.
Denote $[\epsilon_1,\epsilon_2,\ldots,\epsilon_{k-1}]$ by $W$. Applying Lemma 2.3, we have $\sigma(W)=W$.
Suppose $\sigma([\epsilon_k])=[\alpha]$ and $\alpha=\sum_{i=1}^na_i\epsilon_i$.  Since $[\epsilon_k]$ and $W$ are not adjacent, so are  their images under $\sigma$, which implies that $ \alpha\notin W$ and thus
$a_ka_{k+1}\ldots a_n\not=0$.   If $a_{k}\not=0$, then we take $$A =I-a_{k}^{-1}\sum_{i\not=k} a_{ i}E_{ik},$$
 and thus $\theta_{A }\circ\sigma$ fixes   $[\epsilon_k]$.
 If $a_{k}=0$ and $a_{ m}\not=0$ for some $m>k$, let $P_{km}$ be the permutation matrix obtained from $I$ by permuting the $m$-th row and the $k$-th row of $I$, and
   take $$A =(I-a_{ m}^{-1}\sum_{i\notin\{k,m\}} a_{i}E_{ik})\cdot P_{km},$$
   then  $\theta_{A }\circ\sigma$  fixes   $[\epsilon_k]$. Note that $\theta_{A }\circ\sigma$ also fixes $[\epsilon_i]$ for $i=1,2, \ldots, k-1$.\hfill$\square$

 Now, we are ready to give a proof for Theorem 1.3.
 \vskip 2mm
 \noindent{\bf Proof of Theorem 1.3}
 \vskip 1.5mm
  The sufficiency  of Theorem 1.3 is obvious since the product of some automorphisms of a graph is also an automorphism. For the necessity, we complete the proof by establishing several claims. Let $\epsilon_1,\epsilon_2,\ldots, \epsilon_n$ be a set of base of $\mathbb{V}$, $\sigma$   an automorphism of $\mathcal{I}n(\mathbb{V})$.

\vskip 1.5mm
\noindent {\sf Claim 1.}  {\it Either $\sigma$ or $\tau\circ \sigma$ sends $[\epsilon_1]$ to an $1$-dimensional subspace of $\mathbb{V}$.}
\vskip 1.5mm
Clearly, the degree of $\sigma([\epsilon_1])$ must  equal   the degree of $[\epsilon_1]$. Thus by Lemma 2.1, the dimension of $\sigma([\epsilon_1])$ is either $1$ or $n-1$.
If $\sigma$ sends $[\epsilon_1]$ to an $1$-dimensional subspace of $\mathbb{V}$, then there is nothing to do. Otherwise, if $\sigma$ sends $[\epsilon_1]$ to an $(n-1)$-dimensional subspace of $\mathbb{V}$, then  $\tau\circ\sigma$  sends  $[\epsilon_1]$ to $(\sigma([\epsilon_1]))^\perp$, which is an $1$-dimensional subspace of $\mathbb{V}$.

If  $\sigma$  sends  $[\epsilon_1]$ to an $(n-1)$-dimensional subspace of $\mathbb{V}$, we will denote $\tau\circ \sigma$ by $\sigma_1$, otherwise we set $\sigma_1=\sigma$.  Then $ \sigma_1$  sends  $[\epsilon_1]$ to an $1$-dimensional subspace.

\vskip 1.5mm
\noindent {\sf Claim 2.}
There exists an invertible matrix $A $ such that  $\theta_{A }\circ\sigma_1$  fixes  every  $[\epsilon_i]$ for $i=1,2, \ldots, n$.
\vskip 1.5mm
Applying Lemma 2.4, there is an invertible matrix $A_1$ such that $\theta_{A_1}\circ \sigma_1$ fixes $[\epsilon_1]$. Then we can further find an invertible matrix $A_2$ such that $\theta_{A_2}\circ \theta_{A_1}\circ \sigma_1$ respectively fixes $[\epsilon_1]$ and $[\epsilon_2]$. Proceeding in this way, we can find invertible matrices $A_1, A_2, \ldots, A_n$, in sequence, such that  $ \theta_{A_n}\circ\theta_{A_{n-1}}\circ\ldots\circ\theta_{A_1}\circ \sigma_1$   fixes every  $\epsilon_i$ for $i=1,2, \ldots, n$.
Let $A=A_nA_{n-1}\ldots A_1$. Then $\theta_{A }\circ\sigma_1$ is as required.

In the following, we denote $ \theta_{A }\circ\sigma_1$ by $\sigma_2$.
Keep in mind that $\sigma_2$ sends every nontrivial proper subspace to a subspace of equal dimension (by Lemma 2.2)

\vskip 1.5mm
\noindent {\sf Claim 3.} For any $\alpha=\sum_{i=1}^na_i\epsilon_i\not=0$, suppose $\sigma_2([\alpha])=[\beta]$ with $\beta=\sum_{i=1}^nb_i\epsilon_i$. Then $b_k=0$ if and only if $a_k=0$ for $1\leq k\leq n$.
\vskip 1.5mm
 Let $U_k$ denote the subspace spanned by $\epsilon_1, \ldots,  \epsilon_{k-1},\epsilon_{k+1}, \ldots, \epsilon_n$. Then $\sigma_2$ fixes every $U_k$ for $k=1,2, \ldots, n$ (thanks to Lemma 2.3). If $a_k=0$, then $[\alpha]\subset U_k$, which implies that $[\beta]\subset U_k$ and thus $b_k=0$. By considering $\sigma_2^{-1}$ we have $a_k=0$ whenever $b_k=0$.

For $1\leq i<j\leq n$ and $a\in F_q$, Claim 3 shows that $\sigma_2$ sends $[\epsilon_i+a\epsilon_j]$ to a vertex of the form $[\epsilon_i+b\epsilon_j]$ with $b\in F_q$. Thus we can define a function $f_{ij}$ on $F_q$ such that
$f_{ij}(0)=0$ and  $\sigma_2([\epsilon_i+a\epsilon_j])=[\epsilon_i+f_{ij}(a)\epsilon_j]$. Next, we will study the properties of $f_{ij}$.

\vskip 1.5mm
\noindent {\sf Claim 4.} $\sigma_2([\epsilon_i+\sum_{j=i+1}^na_j\epsilon_j])=[\epsilon_i+\sum_{j=i+1}^nf_{ij}(a_j)\epsilon_j]$ for $1\leq i<j\leq n$.
\vskip 1.5mm
Suppose  $\sigma_2([\epsilon_i+\sum_{j=i+1}^na_j\epsilon_j]=[\epsilon_i+\sum_{j=i+1}^nb_j\epsilon_j]$. Applying $\sigma_2$ on
$$[\epsilon_i+\sum_{j=i+1}^na_j\epsilon_j]\subset [\epsilon_i+a_k\epsilon_k, \epsilon_{i+1}, \ldots, \epsilon_{k-1}, \epsilon_{k+1}, \ldots, \epsilon_n], $$
then by Lemma 2.3 we have
$$[\epsilon_i+\sum_{j=i+1}^nb_j\epsilon_j] \subset  [\epsilon_i+f_{ik}(a_k)\epsilon_k, \epsilon_{i+1}, \ldots, \epsilon_{k-1}, \epsilon_{k+1}, \ldots, \epsilon_n],$$
which implies that $b_k=f_{ik}(a_k)$ for any $1\leq i<k\leq n$.
\vskip 1.5mm
\noindent {\sf Claim 5.} Let $2\leq i<j\leq n$. Then \\
(i) $f_{1j}(ab)=f_{1i}(a)f_{ij}(b)$ for any  $a,b\in F_q$;\\
(ii)  $f_{1j}(a)=f_{1i}(a)f_{ij}(1)=f_{1i}(1)f_{ij}(a)$ for any  $a\in F_q$. In particular, $f_{1j}(1)=f_{1i}(1)f_{ij}(1)$. \\
(iii) $\frac{f_{1j}(a)}{f_{1j}(1)} =\frac{f_{1i}(a)}{f_{1i}(1)} =\frac{f_{ij}(a)}{f_{ij}(1)}$.
\vskip 1.5mm
As $[\epsilon_1]\subset [\epsilon_1+a\epsilon_i+ab\epsilon_j, \epsilon_i+b\epsilon_j]$, applying Lemma 2.3 and Claim 4, we have
$$[\epsilon_1]\subset [\epsilon_1+f_{1i}(a)\epsilon_i+f_{1j}(ab)\epsilon_j, \epsilon_i+f_{ij}(b)\epsilon_j],$$
from which it follows that
$$f_{1j}(ab)=f_{1i}(a)f_{ij}(b).$$ Taking $b=1$, we have $$f_{1j}(a)=f_{1i}(a)f_{ij}(1).$$ Similarly, we have $$f_{1j}(b)=f_{1i}(1)f_{ij}(b).$$ If we take $a=b=1$, we have $$f_{1j}(1)=f_{1i}(1)f_{ij}(1),$$  which completes the proof of (ii). (iii)  follows from (ii) immediately.

\vskip 1.5mm
\noindent {\sf Claim 6.} Let $f=\frac{f_{12}}{f_{12}(1)}$. Then\\
(i) $f(a)=\frac{f_{ij}(a)}{f_{ij}(1)}$ for any $a\in F_q$ and $1\leq i<j\leq n$. \\
(ii) $f(ab)=f(a)f(b)$ for any $a,b\in F_q$.\\
(iii) $f(1)=1$  and $f(a^{-1})=f(a)^{-1}$ for $0\not=a\in F_q$.\\
(iv) $f(-1)=-1$ and $f(-a)=-f(a)$ for any $a\in F_q$.\\
(v) $f(a+b)=f(a)+f(b)$ for any $a,b\in F_q$.\\
(vi) $f$ is an automorphism of $F_q$.
\vskip 1.5mm

The definition of $f$ implies that $f(a)=\frac{f_{12}(a)}{f_{12}(1)}$, and (iii) of Claim 5 implies that $$f(a)=\frac{f_{1j}(a)}{f_{1j}(1)}=\frac{f_{ij}(a)}{f_{ij}(1)}, \forall \ 2\leq i<j\leq n,$$ which proves (i) of this claim.

For any $a, b\in F_q$, $$f(ab)=\frac{f_{1n}(ab)}{f_{1n}(1)}=\frac{f_{12}(a)f_{2n}(b)}{f_{12}(1)f_{2n}(1)}=f(a)f(b),$$  which proves (ii) of this claim.

The definition of $f$ implies that $f(1)=1$. If $a\not=0$, from $1=f(1)=f(a^{-1}a)=f(a^{-1})f(a)$ we have $f(a^{-1})=f(a)^{-1}$.

 As $1=f(1)=f((-1)(-1))=f(-1)^2$,  we have $f(-1)=-1$. Further, we have $$f(-a)=f((-1)\cdot a)=f(-1)f(a)=-f(a),$$ which proves (iv).

If $a=0$, (v) is obvious. Assume that $a\not=0$. By $$[\epsilon_1+(a+b)\epsilon_n]\subset [\epsilon_1-a\epsilon_2+a\epsilon_n, \epsilon_2+a^{-1}b\epsilon_n],$$ we have
$$[\epsilon_1+f(a+b)f_{1n}(1)\epsilon_n]\subset [\epsilon_1+f(-a)f_{12}(1)\epsilon_2+f(a)f_{1n}(1)\epsilon_n, \epsilon_2+f(a^{-1}b)f_{2n}(1)\epsilon_n],$$
from which it follows that $$f(a+b)f_{1n}(1)=f(a)f_{1n}(1)+f(a)f_{12}(1)f(a^{-1}b)f_{2n}(1).$$ As $$f_{12}(1)f_{2n}(1)=f_{1n}(1),\ \  f(a)f(a^{-1}b)=f(aa ^{-1}b)=f(b),$$ we further have $f(a+b)=f(a)+f(b)$.

(vi) follows from (ii)-(v) immediately.
\vskip 1.5mm
\noindent {\sf Claim 7.} Set $f_{11}(1)=1$, then $\sigma_2([\epsilon_i+\sum_{j=i+1}^na_j\epsilon_j])=[\epsilon_i+\frac{1}{f_{1i}(1)}\sum_{j=i+1}^nf(a_j)f_{1j}(1)\epsilon_j]$ for $1\leq i<j\leq n$.
\vskip 1.5mm
By Claim 4, $$\sigma_2([\epsilon_i+\sum_{j=i+1}^na_j\epsilon_j])=[\epsilon_i+\sum_{j=i+1}^nf_{ij}(a_j)\epsilon_j], \ \ 1\leq i<j\leq n.$$ Since $f_{ij}(a_{ij})=f(a_{ij})f_{ij}(1)$ and $f_{ij}(1)=\frac{f_{1j}(1)}{f_{1i}(1)}$, the assertion of this claim is confirmed.

\vskip 1.5mm
\noindent {\sf Claim 8.} There is a diagonal matrix $D$ such that $\chi_f^{-1}\circ\theta_D\circ \sigma_2 $ fixes every $1$-dimensional subspace.
\vskip 1.5mm
Let $D=diag(1,f_{12}^{-1}(1),f_{13}^{-1}(1),\ldots, f_{1n}^{-1}(1))$. Then $\theta_D\circ \sigma_2$ sends any  $1$-dimensional subspace $[\epsilon_i+\sum_{j=i+1}^na_j\epsilon_j]$ to $[\epsilon_i+\sum_{j=i+1}^nf(a_j)\epsilon_j]$. Further $\chi_f^{-1}\circ\theta_D\circ \sigma_2 $  fixes every $1$-dimensional subspace. Let $\sigma_3=\chi_f^{-1}\circ\theta_D\circ \sigma_2$.

\vskip 1.5mm
\noindent {\sf Claim 9.} $\sigma_3$ fixes every vertex of $\mathcal{I}n(\mathbb{V})$.
\vskip 1.5mm
Claim 8 has shown that $\sigma_3$ fixes  every $1$-dimensional subspace of $\mathbb{V}$. Now, we proceed by induction on the dimension of subspace to prove that $\sigma_3$ fixes every nontrivial proper subspace of $\mathbb{V}$.
Suppose $\sigma_3$ fixes every $(k-1)$-dimensional subspace of $\mathbb{V}$ and let $W$ be a $k$-dimensional subspace of $\mathbb{V}$ with a basis $w_1,w_2, \ldots, w_k$. The induction hypothesis implies that
$\sigma_3$ fixes  $[w_1,w_2,\ldots, w_{k-1}]$. As $[w_k]\subset W,\  [w_1,w_2,\ldots, w_{k-1}]\subset W$, by applying $\sigma_3$ we have
 $[w_k]\subset \sigma_3(W),\  [w_1,w_2,\ldots, w_{k-1}]\subset \sigma_3(W)$, and thus $W\subseteq \sigma_3(W)$. By comparing their dimensions we have $\sigma_3(W)=W$ .

Claim 9 has  proved that $\sigma_3$ is the identity mapping on  $V$. Thus $\sigma=\tau\circ\theta_X\circ\chi_f$ or $\sigma=\theta_X\circ\chi_f$, where $X=A^{-1}D^{-1}$.

Suppose $$\sigma=\tau^{\delta_1}\circ\theta_{X_1}\circ\chi_{f_1}=\tau^{\delta_2}\circ\theta_{X_2}\circ\chi_{f_2}$$
 are two decompositions  of $\sigma$, where $\delta_i=1$ or $0$, $X_1, X_2$ are $n\times n$ invertible matrices over $F_q$ and $f_1, f_2$ are automorphisms of $F_q$. We first prove that $\delta_1=\delta_2$.
 Indeed, if $\delta_1\not= \delta_2$, say $\delta_1=1$ and $\delta_2=0$, then $\tau^{\delta_1}\circ\theta_{X_1}\circ\chi_{f_1}$ sends $[\epsilon_1]$ to an $(n-1)$-dimensional subspace, however $\tau^{\delta_2}\circ\theta_{X_2}\circ\chi_{f_2}$  sends $[\epsilon_1]$ to an $ 1$-dimensional subspace, a contradiction.
 Thus, $\delta_1=\delta_2$ and $$\theta_{X_1}\circ\chi_{f_1}=\theta_{X_2}\circ\chi_{f_2},$$ which further implies that $$\theta_{X_2^{-1}X_1} = \chi_{f_2f_1^{-1}}.$$
  As $\theta_{X_2^{-1}X_1}$ fixes every $[\epsilon_i]$ for $i=1,2, \ldots,n$, we find that $X_2^{-1}X_1$ is a diagonal matrix. Furthermore, since   $\theta_{X_2^{-1}X_1}$ also fixes every $[\epsilon_i+\epsilon_j]$ for $i\not=j$, we find that the diagonal matrix must be a nonzero scalar matrix. Hence $X_2$ must be a nonzero scalar multiple of $X_1$ and thus $\theta_{X_2}=\theta_{X_1}$, which implies that $\chi_{f_2}=\chi_{f_1}$.\hfill$\square$

Keeping in mind that the automorphism group of $F_q$, with $q=p^m$, is a cyclic group of order $m$, then one   can easily conclude  Corollary 1.4  by applying Theorem 1.3.

\end{document}